\documentclass[12pt]{amsart}

\newif\ifpictures
\picturestrue

\usepackage{amsmath}
\usepackage{amssymb}
\usepackage{latexsym}
\usepackage[dvips]{graphicx}
\ifpictures
  \input{epsf}
\fi

\newtheorem{theorem}{Theorem}

\newtheorem{corollary}[theorem]{Corollary}

\newtheorem{proposition}[theorem]{Proposition}
\newtheorem{conjecture}[theorem]{Conjecture}

\newenvironment{definition}{ \noindent {\bf{Definition.}} }{}

\newcommand{\R}{\mathbb{R}}

\newcommand{\N}{\mathbb{N}}
\renewcommand{\P}{\mathbb{P}}

\copyrightinfo{}{}

\begin{document}


\title{Real $k$-flats tangent to quadrics in $\R^n$}

\author{Frank Sottile}
\address{Department of Mathematics\\
          Texas A\&M University\\
         College Station\\
         Texas \ 77843\\
         USA}
\email{sottile@math.tamu.edu}
\urladdr{http://www.math.tamu.edu/$\sim$sottile}

\author{Thorsten Theobald}
\address{Institut f\"ur Mathematik, MA 6-2, 
  Technische Universit\"at Berlin,
  Stra{\ss}e des 17.~Juni 136, D--10623 Berlin,
  Germany}
\email{theobald@math.tu-berlin.de}
\urladdr{http://www.math.tu-berlin.de}

\thanks{Research of Sottile supported by NSF CAREER grant DMS-0070494 and
  the Clay Mathematical Institute}
\keywords{
tangents, transversals, quadrics, enumerative geometry, real solutions, \\
\indent \hspace{113pt}Grassmannian.
}

\begin{abstract}
 Let $d_{k,n}$ and $\#_{k,n}$ denote the dimension and the degree of
 the Grassmannian $\mathbb{G}_{k,n}$, respectively.
 For each $1 \le k \le n-2$ there are $2^{d_{k,n}} \cdot \#_{k,n}$
 (a priori complex) $k$-planes in $\P^n$
 tangent to $d_{k,n}$ general quadratic hypersurfaces in $\P^n$. 
 We show that this class of enumerative problems is fully real, i.e., 
 for $1 \le k \le n-2$ there
 exists a configuration of $d_{k,n}$ real quadrics in (affine) real space
 $\R^n$ so that all the mutually 
 tangent $k$-flats are real.
\end{abstract}

\maketitle

\section*{Introduction}

Understanding the real solutions of a system of polynomial equations is a
fundamental problem in mathematics (see, e.g., \cite{sturmfels-cbms} 
for some recent lines of research and applications). 
However, as pointed out in \cite[p.~55]{fulton-b96}, 
even for problem classes with a finite number of complex solutions 
(enumerative problems), the question of how many solutions can be 
real is still widely open. 
A class of enumerative problems is called \emph{fully real} if there are general
real instances for which all the (a priori  complex) solutions are
real.

One of us (Sottile) began a systematic study of this question in the
\emph{special Schubert calculus}~\cite{So97a,So99a}, a 
class of enumerative problems from classical algebraic geometry.
This special Schubert calculus asks for linear subspaces
of a fixed dimension meeting some given (general) linear subspaces 
(whose dimensions and number ensure a finite number of solutions)
in $n$-dimensional complex projective space $\mathbb{P}^n$.
For any given dimensions of the subspaces, this problem is fully real,
i.e., there exist \emph{real} linear subspaces for which each of the
a priori complex solutions is \emph{real}.
In particular, for $1 \le k \le n-2$ there are 
$d_{k,n}:=(k+1)(n-k)$ real $(n{-}k{-}1)$-planes 
$U_{1}, \ldots, U_{d_{k,n}}$ in $\mathbb{P}^n$ with
\[
  \#_{k,n} := \frac{1!2!\cdots k!((k+1)(n-k))!}{(n-k)!(n-k+1)! \cdots n!}
\]
real $k$-planes meeting $U_1, \ldots, U_{d_{k,n}}$. 
Here, $d_{k,n}$ and $\#_{k,n}$ are the dimension and the degree of the
Grassmannian $\mathbb{G}_{k,n}$, respectively
(see~\cite{kleiman-laksov-72,schubert-1886}). 
These were the first results showing that a large class of 
non-trivial enumerative problems is fully real.
Recently, Vakil~\cite{Va04} has shown that any Schubert problem on
a Grassmannian is fully real.

We continue this line of research by considering $k$-flats tangent to quadratic
hypersurfaces (hereafter \emph{quadrics}).  
This is also motivated by recent 
investigations in computational geometry (see~\cite{MPT01,M2,STh02}).
It was shown in~\cite{STh02} that $2n{-}2$ general \emph{spheres} 
in affine real space $\mathbb{R}^n$ have at most $3\cdot 2^{n-1}$ 
common tangent lines in $\mathbb{C}^n$, and that 
there exist spheres for which all the
a priori complex tangent lines are \emph{real}.
The present paper addresses the following question:
What is the maximum number of real $k$-flats simultaneously tangent to
$d_{k,n}$ general quadrics in $\mathbb{R}^n$ 
(respectively in $\mathbb{P}_{\mathbb{R}}^n$)? 
As this problem may be formulated as the complete intersection of 
$d_{k,n}$ quadratic equations on the Grassmannian of $k$-planes in $\P^n$,
the expected number of complex solutions is the product of the degrees of
the equations with the degree of the Grassmannian, i.e.,
$ 2^{d_{k,n}}\cdot \#_{k,n}$. We show that the problem is fully real:

\begin{theorem}\label{T:allreal}
 Let $1 \le k \le n-2$.
 Given $d_{k,n}$ general quadrics in $\P^n$ there are 
 $2^{d_{k,n}} \cdot \#_{k,n}$ complex $k$-planes 
 that are simultaneously tangent to all $d_{k,n}$ 
 quadrics. Furthermore, there is a choice of quadrics
 in $\R^n$ for which all the
 $k$-flats are real, distinct, and lie in affine space $\R^n$.
\end{theorem}

When $k=1$, we have
$d_{1,n} = 2(n-1)$ and $\#_{1,n}$ is the Catalan number
$\#_{1,n} = \frac{1}{n} \binom{2n-2}{n-1}$.
The following table exhibits 
the large discrepancy between the number of 
lines tangent to spheres and the number of lines tangent to general quadrics.
When $n=3$ this discrepancy was accounted for by Aluffi and Fulton~\cite{AF}.
\[
  \begin{tabular}{|c||c|c|c|c|c|c|c|}\hline
    $n$& 3 & 4 & 5 & 6 & 7 & 8 & 9\\\hline\hline \rule{0pt}{13pt}
    $3\cdot 2^{n-1}$&12&24&48&96&192&384&768\\\hline\rule{0pt}{12pt}
    $2^{d_{1,n}} \cdot \#_{1,n}$&32&320&3584&43008&540672&7028736&93716480\\\hline
  \end{tabular}
\]

In Section~\ref{S:plucker}, we review some facts on Pl\"ucker coordinates
of $k$-planes in projective space.
In Section~\ref{S:proof}, we combine recent results in the real 
Schubert calculus with classical perturbation arguments adapted 
to the real numbers to prove Theorem~\ref{T:allreal}.
Since the proof for general $(k,n)$ is non-constructive, we give a 
symbolic, constructive
proof for the case $(k,n) = (1,3)$ in Section~\ref{S:symbolic}.

\section{Preliminaries}\label{S:plucker}

We review the well-known \emph{Pl\"ucker coordinates} of $k$-dimensional
linear subspaces (hereafter $k$-planes) in complex projective space
$\P^n$ (see, e.g., \cite{hodge-pedoe}).
Let $U$ be a $k$-plane in $\P^n$ which is spanned by the columns of
an $(n+1) \times (k+1)$-matrix $L$.
For every subset $I \subset \{0, \ldots, n\}$ of size $k+1$ let
$p_{I}$ be the $(k+1) \times (k+1)$-subdeterminant 
of $L$ given by the rows in $I$
and let $N:=\binom{n+1}{k+1}-1$.
Then $p := (p_I)_{I \subset \{0, \ldots, n\}, |I| = k+1} \in \P^N$
is the \emph{Pl\"ucker coordinate} of $U$.
The set of all $k$-planes in $\P^n$ is called the \emph{Grassmannian
of $k$-planes in} $\P^n$ and is denoted by $\mathbb{G}_{k,n}$.
If the indices are written as ordered tuples then the Pl\"ucker
coordinates are skew-symmetric in the indices.
$\mathbb{G}_{k,n}$ is in 1-1-correspondence with the set of
vectors in $\P^N$ satisfying the \emph{Pl\"ucker relations}, i.e.,
\begin{equation}
  \label{eq:plueckerrelation}
  \sum_{l=1}^{k+1} (-1)^l p_{i_1 \ldots \hat{i_l} \ldots i_{k+1}} \,
     p_{j_1 \ldots j_{k-1} i_l} \ =\ 0
\end{equation}
for every $I=\{i_1, \ldots, i_{k+1}\}$, 
$J = \{j_1 , \ldots, j_{k-1}\} \subset \{0, \ldots, n\}$
of strictly ordered index sets
(where $\hat{\hspace*{1mm}}$ over an index means that it is omitted).
See, e.g., \cite[\S VII.6]{hodge-pedoe}.
By Schubert's results \cite{schubert-1886}, the dimension of 
$\mathbb{G}_{k,n}$ is $d_{k,n} = (k+1)(n-k)$ and its degree is $\#_{k,n}$.

If an $(n{-}k{-}1)$-plane $V$ is given as the intersection of the
$k+1$ hyperplanes
$\sum_{i=0}^n v_i^{(0)} x_i {~=~} 0, \ldots , \sum_{i=0}^n v_i^{(k)} x_i {~=~} 0$, 
then the \emph{dual Pl\"ucker
coordinate} $q = (q_I)_{I \subset \{0, \ldots, n\}, |I| = k+1} \in \P^N$ of 
$V$ is defined by the $(k+1) \times (k+1)$-subdeterminants of the matrix
with columns $v^{(0)}, \ldots, v^{(k)}$.

A $k$-plane $U$ intersects an $(n{-}k{-}1)$-plane $V$ in $\P^n$
if and only if the dot product of the Pl\"ucker coordinate $p$ of $U$ and 
the dual Pl\"ucker coordinate $q$ of $V$ vanishes, i.e., if and only if
\begin{equation}
  \label{eq:transversalcond}
  \sum_{I \subset \{0, \ldots, n\}, |I| = k+1} p_{I} q_{I} \ = \ 0
\end{equation}
(see, e.g., \cite[Theorem~VII.5.I]{hodge-pedoe}).

We use Pl\"ucker coordinates to characterize the $k$-planes tangent to a 
given quadric in $\mathbb{P}^n$ 
(see \cite{M2}). 
We identify a quadric $x^{\mathrm{T}} Q x = 0$ in $\P^n$ with its
symmetric  $(n{+}1)\times(n{+}1)$-representation matrix $Q$. 
Further, for $r \in \N$
let $\wedge^r$ denote the $r$-th exterior power of matrices 
\begin{eqnarray*}
  \wedge^r \ \colon\ \mathbb{C}^{m\times n} & \to & 
             \mathbb{C}^{\binom{m}{2}\times\binom{n}{2}}
\end{eqnarray*}
(see~\cite{M2}).
The row and column indices of the resulting matrix are subsets of 
cardinality~$r$ of $\{1, \ldots, m\}$ and $\{1, \ldots, n\}$, respectively.
For $I \subset \{1, \ldots, m\}$ with $|I| = r$ and
$J \subset \{1, \ldots, n\}$ with $|J| = r$,
$\bigl(\wedge^r A\bigr)_{I,J}$ is the subdeterminant of $A$ whose
rows are indexed by $I$ and whose columns are indexed by $J$.
If a $k$-plane $U \subset \P^n$ is spanned by the columns of 
an $(n{+}1)\times (k+1)$-matrix $L$,
then the $\binom{n+1}{k+1} \times 1$-matrix $\wedge^{k+1} L$,
considered as a vector in $\P^N$, 
is the Pl\"ucker coordinate of $U$. 

Recall the following algebraic characterization of tangency:
A $k$-plane $U$ is tangent to a quadric $Q$ if the
restriction of the quadratic form to $U$ is singular
(which includes the case $U \subset Q$).
When the quadric is smooth, this implies that $U$ is tangent to the
quadric in the usual geometric sense.

\begin{proposition} [Proposition 5.5.3 of~\cite{M2}]
\label{le:quadrictangentcond}
 A $k$-plane $U \subset \P^n$ is tangent to a quadric $Q$ if and only
 if the Pl\"ucker coordinate $p_U$ of $U$ satisfies
\begin{equation}
  \label{eq:quadrictangentcond}
  p_U^{\mathrm{T}} \, \bigl(\wedge^{k+1} Q\bigr) \, p_U = 0 \, .
\end{equation}
\end{proposition}

A $k$-\emph{flat} in affine real space 
$\R^n$ is a $k$-dimensional affine subspace in $\R^n$.
Throughout the paper we assume that $\R^n$ is naturally embedded in 
the real projective space $\P_{\mathbb{R}}^n$ via 
$(x_1, \ldots, x_n) \mapsto (1,x_1, \ldots, x_n) \in \P_{\mathbb{R}}^n$.

\section{Proof of the main theorem}\label{S:proof}

We first illustrate the essential geometric idea underlying our
constructions for $(k,n) = (1,3)$, which is the first nontrivial case.
Here, Theorem~\ref{T:allreal} 
states that there exists a configuration of four quadrics in $\R^3$ with $32$
distinct real common tangent lines. 

By~\eqref{eq:transversalcond},
the set of lines meeting four given lines in $\P^3$ is the intersection
of four hyperplanes on the Grassmannian $\mathbb{G}_{1,3}$, and hence
there are at most two
or infinitely many common lines meeting $\ell_1, \ldots, \ell_4$.
If $e_1$ and $e_2$ are opposite edges in a tetrahedron $\Delta \subset \R^3$,
then the lines underlying $e_1$ and $e_2$ are the two common transversals 
of the four lines underlying the other four edges 
(see~Figure~\ref{fi:4lines2transversals}).

\ifpictures
\begin{figure}[htb]
\[
  \unitlength1.2pt
  \begin{array}{c@{\hspace*{1.5cm}}c}
  \begin{picture}(150,130)(0,12)
   \put(3,3){\epsfxsize=168pt\epsfbox{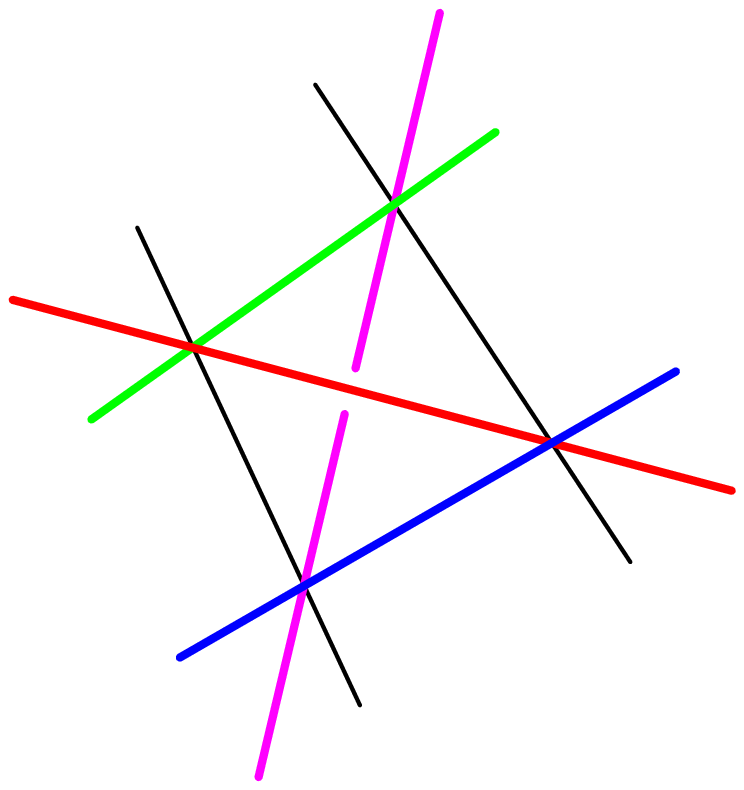}}
   \put( 41,  7){$\ell_1$}
   \put( 85, 40){$\ell_2$}
   \put(  8, 96){$\ell_3$}
   \put( 50,105){$\ell_4$}
   \put( 35, 55){$e_1$}
   \put( 98, 93){$e_2$}
  \end{picture}
& \\ [-5.8cm] &
  \epsfxsize=2.4in \epsfbox{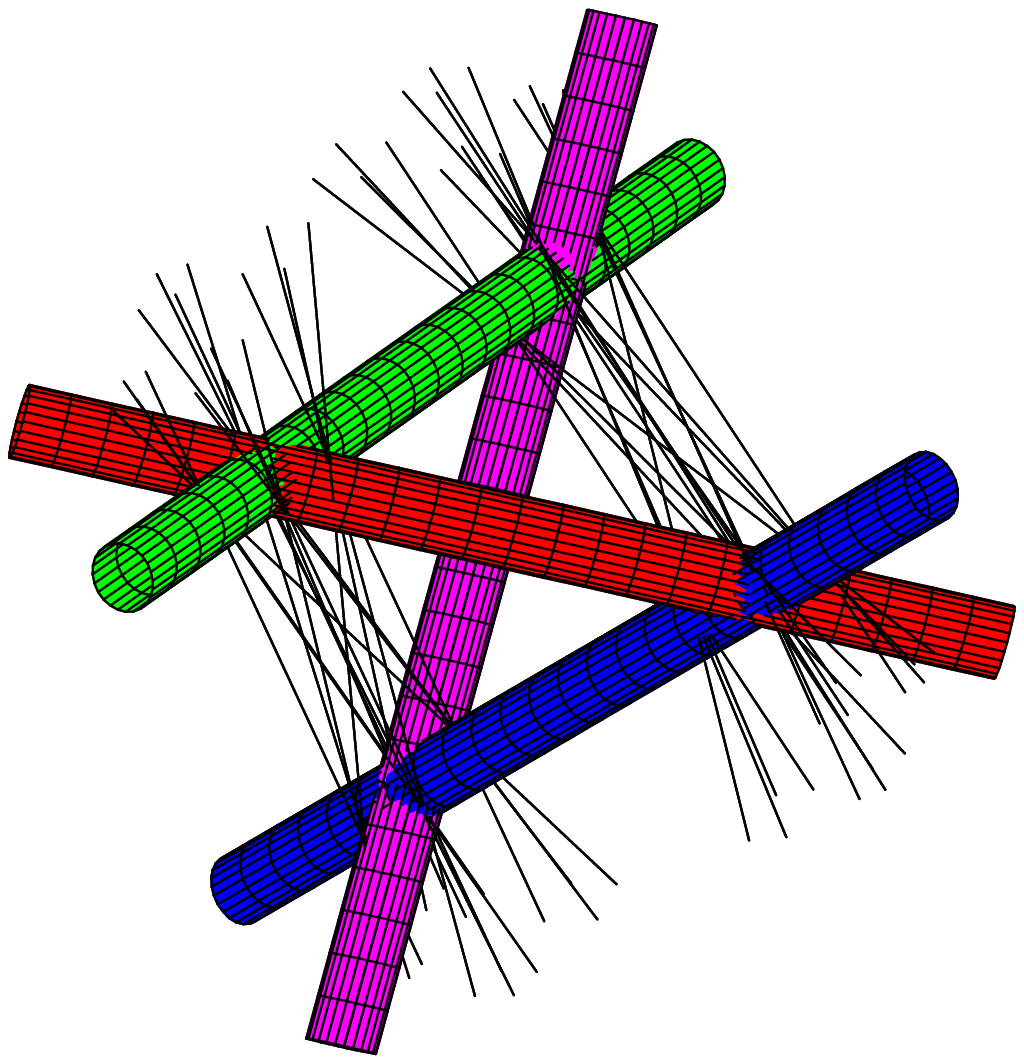} \\
\end{array}
\]
\caption{A tetrahedron configuration of four lines in $\R^3$ with two 
real transversals and a configuration of four quadrics with 32 real tangents.}
\label{fi:4lines2transversals}
\end{figure}
\fi

Consider the lines $\ell_1, \ldots, \ell_4$ as (degenerate)
infinite circular cylinders with radius $r=0$. 
When the radius is slightly increased, then the cylinders 
intersect pairwise in the regions (combinatorially)
given by the four vertices of $\Delta$, and the common tangents roughly 
have the direction of $e_1$ or $e_2$. 
Since the neighborhood of a vertex is divided into four regions 
by the two cylinders, and since each region contains common tangents,
this gives $4 \cdot 4$ tangents close to the direction of
$e_1$ and $4 \cdot 4$ tangents close to the direction of $e_2$
(see Figure~\ref{fi:4lines2transversals}.)

For the general case, let $1 \le k \le n-2$. By Section~\ref{S:plucker},
the number of $k$-planes in $\P^n$ simultaneously 
meeting $d_{k,n}$ general $(n{-}k{-}1)$-planes is $\#_{k,n}$.
We begin with a configuration of $d_{k,n}$ real $(n{-}k{-}1)$-flats 
$U_1, \ldots, U_{d_{k,n}}$ in $\R^n$ having 
$\#_{k,n}$ real $(n{-}k{-}1)$-flats simultaneously meeting
$U_1, \ldots, U_{d_{k,n}}$.
We then argue that we can replace each of these $(n{-}k{-}1)$-flats by a real
quadric such that for each of the $k$-flats,
there are $2^{d_{k,n}}$ nearby real $k$-flats tangent to each quadric.

\begin{proposition}
\label{pr:realflats}
 For $1 \le k \le n-2$, there exists a configuration of $d_{k,n}$ real 
$(n{-}k{-}1)$-flats $U_1, \ldots, U_{d_{k,n}}$ in 
$\R^n$ such that there exist exactly $\#_{k,n}$ real $k$-flats
simultaneously meeting $U_1, \ldots, U_{d_{k,n}}$.
\end{proposition}

\begin{proof}
 The corresponding statement for real projective space $\P_{\R}^n$ was proven
 for $k=1$ in~\cite[Theorem~C]{So97a} and for $k\ge 2$ in~\cite{So99a}. 
 We deduce the affine counterpart above simply by removing a real hyperplane
 that contains none of the $(n{-}k{-}1)$-flats or any of the transversal
 $k$-flats. 
\end{proof}

For $k=1$, the purely existential statement in~\cite{So97a} 
and Proposition~\ref{pr:realflats}
was improved
 by Eremenko and Gabrielov~\cite{EG02} who gave the following explicit
 construction of such a collection of $(n{-}2)$-flats.
 Let $\gamma:\R \to \R^n$, 
$\gamma(s)\ :=\ (1,s,s^2,\dotsc,s^{n-1})^{\mathrm{T}}$
 be the moment curve in $\R^n$.
 For each $s\in\R$, set $U(s)$ to be 
\[
  U(s)\ :=\ \textrm{affine span}\bigl(\gamma(s),\gamma'(s),\dotsc,
                \gamma^{(n{-}3)}(s)\bigr)\,.
\]
 Geometrically, $U(s)$ is the $(n{-}2)$-flat osculating 
 the moment curve at the point $\gamma(s)$.
 By \cite{EG02}, for any distinct
 $s_1,\ldots,s_{2n-2}\in\R$, the $(n{-}2)$-flats 
 $U(s_1), U(s_2),\dotsc, U(s_{2n-2})$ have 
 exactly $\#_{1,n} = C_{n-1}$ common real transversals, where
 $C_{n} = \frac{1}{n+1} \binom{2n}{n}$ is the $n$-th Catalan number.
For general $k$, it is only known that 
there \emph{exist} distinct $s_1, \ldots, s_{d_{k,n}} \in \R$
such that there are $\#_{k,n}$ distinct real $k$-flats meeting 
the osculating $(n{-}k{-}1)$-flats to the moment curve
at $s_1, \ldots, s_{d_{k,n}}$ \cite{So99a}.
The \emph{conjecture on total reality} in 
\cite[\S 1 and \S 4]{sedykh-shapiro-2002}
conjectures that any choice of distinct 
$s_1, \ldots, s_{d_{k,n}} \in \R$ implies reality of all transversal subspaces.

\medskip

\begin{definition}
{\rm Suppose that $1 \le k\le n-2$, and let $U \subset \R^n$ be a $k$-flat and 
 $r > 0$.
 The {\it $k$-cylinder} $\textrm{Cy}(U,r)$ is 
 the set of points having Euclidean distance $r$ from $U$.
}
\end{definition}

\medskip 

This quadratic hypersurface is smooth in $\R^n$ but its extension
to $\P^n$ is singular.
A $k'$-flat $V \subset \R^n$ is tangent to $\textrm{Cy}(U,r)$ 
if and only if its Euclidean distance to $U$ is $r$.

We will use the following basic property of intersection
multiplicities~\cite[p.~1]{fulton-b96}.

\begin{proposition}\label{pr:multiplicity}
 Let $\mathcal{A}$ be an algebraic curve in $\P^n$, 
 and let $x$ be a singular point on $\mathcal{A}$. 
 For any hyperplane $H \subset \P^n$ 
 such that $x$ is an isolated point in $\mathcal{A} \cap H$,
 the intersection multiplicity of $\mathcal{A}$ 
 and $H$ in $x$ is greater than $1$.
\end{proposition}

\begin{theorem}
 \label{th:inductive}
 Let $1 \le k \le n-2$, and let
 $U_1,U_2,\ldots,U_{d_{k,n}}$ be $(n{-}k{-1})$-flats in
 $\R^n$ having exactly $\#_{k,n}$ common transversal $k$-flats, all real.
 For each $i=0,1,\dotsc, d_{k,n}$, there 
 exist $r_1,\dotsc,r_{i} > 0$ such that there
 are exactly $2^i \cdot \#_{k,n}\,$ distinct $k$-flats, each of them real, 
 that are simultaneously tangent to each of the 
 $(n{-}k{-}1)$-cylinders ${\rm{Cy}}(U_j,r_j)$,
 $j=1,\dotsc,i$, and also meet each of the $(n{-}k{-}1)$-flats
 $U_{i+1},\dotsc,U_{d_{k,n}}$. 
\end{theorem}

The case of $i=d_{k,n}$ implies Theorem~\ref{T:allreal}.

\begin{proof}
  We induct on $i$, with the case of $i=0$ being the hypothesis of the theorem.

  Suppose that $i\leq d_{k,n}$ and that there exist $r_1,\dotsc,r_{i-1} > 0$
  such that there are exactly $2^{i-1} \cdot \#_{k,n} \,$ distinct $k$-flats 
  $V_1, \ldots, V_{2^{i-1} \#_{k,n}}$ which 
  are simultaneously tangent to $\mathrm{Cy}(U_j, r_j)$
  for each $j=1, \ldots, i-1$,
  and meet each of $U_{i}, \ldots, U_{d_{k,n}}$, and each of these
  $k$-flats is real.

  Now we drop the condition that the $k$-flats meet $U_i$.
  Let $\mathcal{A} \subset \mathbb{G}_{k,n}$ be the curve 
  of $k$-flats that are tangent to the cylinders
  $\textrm{Cy}(U_j,r_j)$ for $j=1,\dotsc,i{-}1$ and that also 
  meet each of the $(n{-}k{-}1)$-flats $U_{i+1},\dotsc,U_{d_{k,n}}$.
  Since $\mathcal{A}$ is the intersection of $i{-}1$ quadrics (the tangency conditions)
  with $d_{k,n}-i$ hyperplanes (conditions to meet the remaining $U_j$) on
  the Grassmannian, it has degree at most $2^{i-1} \#_{k,n}$.
  Since its intersection with the hyperplane defined by $U_i$ consists
  of $2^{i-1} \#_{k,n}$ points, we conclude that the degree of $\mathcal{A}$ is 
  $2^{i-1} \#_{k,n}$ and (by Proposition~\ref{pr:multiplicity}) that
  each of these points is a smooth point of $\mathcal{A}$.

  Let $V \in \{V_1, \ldots, V_{2^{i-1} \#_{k,n} }\}$.
  Since $V$ is a smooth real point of the real curve 
  $\mathcal{A} \subset \mathbb{G}_{k,n}$ (i.e., 
  $V \subset \P_{\R}^n$),
  the real points of $\mathcal{A}$ contain a smooth arc $\alpha$
  containing $V$ with 
  $\alpha \cap (\{V_1,\dotsc,V_{2^{i-1} \#_{k,n}}\} \setminus V) = \emptyset$.
  Let $\varphi\colon(-\delta,\delta)\to\alpha$ be a smooth parametrization of
  the arc $\alpha$ with $\varphi(0)=V$.
  Such a parametrization exists, for example, by the Implicit Function Theorem.

  Thus, for $t\in(-\delta,\delta)\setminus \{0\}$, 
  the real $k$-flat 
  $\varphi(t)$ does not meet $U_i$ and so it has a positive
  Euclidean distance $d(t)$ from $U_i$.
  Since $d(t)$ is a continuous function of $t$, 
  for $\rho \in \R$ with $0 < \rho < \min\{d(-\delta/2),d(\delta/2)\}$ 
  there are at least two distinct real $k$-flats in $\alpha$
  whose Euclidean distance to
  $U_i$ is $\rho$.

  In this way, we obtain $2^{i-1} \cdot \#_{k,n}$ such arcs, each containing 
   one of $V_1,\dotsc,V_{2^{i-1} \#_{k,n}}$.
  We may assume that these arcs are pairwise disjoint.
  Let $0<r_i$ be small enough to ensure that each arc contains
  two $k$-flats 
  having Euclidean distance $r_i$ from $U_i$.
  This gives {\it at least} $2 \cdot 2^{i-1} \cdot \#_{k,n}$ real $k$-planes in 
  $\mathcal{A}$ whose
  Euclidean distance to $U_i$ is $r_i$. 
  Since $2^i \cdot \#_{k,n}$ is the maximum number of $k$-flats with this property,
  there are exactly $2^i \cdot \#_{k,n}$ distinct 
  $k$-flats tangent to  $\textrm{Cy}(U_j,r_j)$ for 
  $j=1,\dotsc,i$ and that also meet each of the $(n{-}k{-}1)$-flats 
  $U_{i+1},\dotsc,U_{d_{k,n}}$.  
\end{proof} 

Since the number of
real $k$-flats will not change under a small perturbation of the 
$k$-cylinders
$\textrm{Cy}(U_j,r_j)$, we may replace them by quadrics
which are smooth in $\P^n$. Let $\text{sign}(Q)$ 
denote the signature of a quadric $Q \subset \P^n$.

\begin{corollary}
\label{co:realitycol}
Let $1 \le k \le n-2$. 
For 
\[
 (s_1, \ldots, s_{d_{k,n}}) \in 
  \begin{cases}
    \{n-1, n-3, \ldots, 2k-n+1\}^{d_{k,n}} & \text{if $\; k \ge n/2$}\, , \\
    \{n-1, n-3, \ldots, 2 \cdot 
     \left(\frac{n-1}{2} - \lfloor \frac{n-1}{2} \rfloor \right)\}^{d_{k,n}} & \text{if $\; k < n/2$} \, 
  \end{cases}
\]
there exist smooth quadrics $Q_1, \ldots, Q_{d_{k,n}} \subset \P_{\R}^n$ 
with $|\mathrm{sign}(Q_i)| = s_i$, $1 \le i \le d_{k,n}$,
such that
the $\#_{k,n}$ (complex) common tangent $k$-flats to $Q_1, \ldots, Q_{d_{k,n}}$
are all real, distinct, and lie in affine space $\R^n$.
\end{corollary}

\begin{proof}
Since the absolute value of the 
signature of an $(n-k-1)$-cylinder is $k$, the proof
immediately follows from the possible perturbations of the quadratic
form in $\P^n$ of the type
\[
  -r^2 x_0^2 + x_1^2 + \dotsb + x_{k+1}^2\,. 
\] 
\vspace*{-7ex}

\end{proof}

We conjecture that the reality statement holds 
for signatures not covered by Corollary~\ref{co:realitycol}.

\begin{conjecture}
Let $1 \le k \le n-2$. For
\[
\begin{array}{c}
(s_1, \ldots, s_{d_{k,n}}) \in
\{n-1, n-3, \ldots, 2 \cdot 
     \left(\frac{n-1}{2} - \lfloor \frac{n-1}{2} \rfloor \right)\}^{d_{k,n}}
\end{array}
\]
there exist smooth quadrics $Q_1, \ldots, $ $Q_{d_{k,n}} \subset \P_{\R}^n$ 
with $|\mathrm{sign}(Q_i)| = s_i$, $1 \le i \le d_{k,n}$,
such that
the $\#_{k,n}$ (complex) common tangent $k$-flats to $Q_1, \ldots, Q_{d_{k,n}}$
are all real, distinct, and lie in affine space $\R^n$.
\end{conjecture}

The first case of this conjecture which is not covered by
Corollary~\ref{co:realitycol} is when $k=3$ and $n=5$ and the signature is zero.
That is, for $3$-flats tangent to 8 smooth quadrics in $\P_{\R}^5$, 
with at least one having
signature zero.
We remark that an argument perturbing cylinders to singular quadrics gives an
analog to Corollary~\ref{co:realitycol} concerning $k$-flats tangent to singular
quadrics. 
We omit its complicated formulation.

\section{A constructive proof for lines in dimension~3} \label{S:symbolic}

Our proof of Theorem~\ref{T:allreal} was non-constructive. 
We close this paper by providing a 
constructive proof in the first nontrivial case,
$(k,n) = (1,3)$, i.e., the
real lines tangent to four quadrics in 3-space.
In order to realize the tetrahedral configuration of
Figure~\ref{fi:4lines2transversals} in $\P_{\R}^3$,
let $\ell_1, \ldots, \ell_4$ be given by the following
equations:
\begin{equation}\label{E:proj-lines}
 \begin{array}{rcll}
  \ell_1\,:\,x_0=x_3=0 \, ; \quad
  \ell_2\,:\,x_0=x_1=0 \, ; \quad
  \ell_3\,:\,x_1=x_2=0 \, ; \quad
  \ell_4\,:\,x_2=x_3=0 \, .
 \end{array}
\end{equation}
The two common transversal lines are given by $x_2=x_4 = 0$
and by $x_1 = x_4 = 0$.

For parameters $\alpha,\beta\in\R$, consider the four quadrics
\[
 \begin{array}{ccrcl}
  Q_1 & : & x_0^2 + x_3^2 - \beta  (x_1^2 + x_2^2) & = & 0\rule{0pt}{12pt}\, , \\
  Q_2 & : & x_0^2 + x_1^2 - \beta  (x_2^2 + x_3^2) & = & 0\rule{0pt}{12pt}\, , \\
  Q_3 & : & x_1^2 + x_2^2 - \alpha (x_0^2 + x_3^2) & = & 0\rule{0pt}{12pt}\, , \\
  Q_4 & : & x_2^2 + x_3^2 - \alpha (x_0^2 + x_1^2) & = & 0\rule{0pt}{12pt}\, .
 \end{array}
\]
For $\alpha=\beta=0$, the four quadrics become the corresponding lines in
$\mathbb{P}_{\R}^3$. 
For small $\alpha,\beta>0$, these quadrics are deformations
of the lines with rank~4 and signature~0---smooth ruled surfaces.

\begin{theorem}
\label{th:32construction} 
 Let $(\alpha,\beta) \in \R^2$ satisfy
\[
 \alpha\beta(1-\alpha\beta)(1-\beta^2)(1-\alpha^2)
  \bigl((1-\alpha)^2(1-\beta)^2-16\alpha\beta\bigr) \ \neq \ 0\, .  
\]
 Then there are $32$ distinct (possibly complex) common tangent lines to
 $Q_1, \ldots, Q_4$.
 If\/ $0<\alpha,\beta<3 - 2 \sqrt{2}$, then each of 
 these $32$ tangent lines is real.
\end{theorem}

\begin{proof}
 Since the quadrics only contain monomials of the form $x_i^2$,
 the tangent equations~\eqref{eq:quadrictangentcond} of $Q_1, \ldots, Q_4$
 only contain monomials of the form $p_{ij}^2$. 
 Hence, the four tangent equations give the following system of
 linear equations in $p_{01}^2, \ldots, p_{23}^2$:
\[
  \left(
  \begin{matrix}
    -\beta   &-\beta  &       1   &  \beta^2 & - \beta & - \beta \\
           1 &-\beta  & -\beta    & -\beta   & - \beta & \beta^2  \\
    -\alpha  &-\alpha &  \alpha^2 &        1  & - \alpha & - \alpha \\
    \alpha^2 &-\alpha & -\alpha   & -\alpha   & - \alpha &        1 
  \end{matrix}
  \right)
  \left(
  \begin{matrix}
    p_{01}^2 \\
    p_{02}^2 \\
    p_{03}^2 \\
    p_{12}^2 \\
    p_{13}^2 \\
    p_{23}^2
  \end{matrix}
  \right)
  \ =\  0 \, .
\]

Permute the variables into the order
$(p_{02},p_{13},p_{03},p_{12},p_{01},p_{23})$.
For $\alpha,\beta \in \R$ satisfying 
 \begin{equation}\label{E:abnez}
   \alpha\beta(1-\alpha\beta)(1+\beta)(1+\alpha)\ \neq \ 0\,,
 \end{equation}
Gaussian elimination yields the following system:
\[
  \left(
  \begin{matrix}
    -\beta & -\beta & (1-\alpha)(1-\beta) & 0  & 0  & 0 \\
    0      & 0      & \alpha & -\beta & 0  & 0 \\
    0      & 0      & 0      & -\beta & \alpha & 0 \\
    0      & 0      & 0      &   0    & \alpha & -\beta
  \end{matrix}
  \right)
  \left(
  \begin{matrix}
    p_{02}^2 \\
    p_{13}^2 \\
    p_{03}^2 \\
    p_{12}^2 \\
    p_{01}^2 \\
    p_{23}^2
  \end{matrix}
  \right)
  \ =\  0 \, .
\]
Together with the Pl\"ucker equation~\eqref{eq:plueckerrelation}, 
this gives the following system of equations: 
 \begin{eqnarray}
    - \beta p_{02}^2 - \beta p_{13}^2
    + (1-\alpha)(1-\beta) p_{03}^2 & = & 0 \, , \label{eq:eq1} \\
    p_{01} p_{23} - p_{02} p_{13} + p_{03} p_{12} & = & 0 \, , \label{eq:eq2} \\
    \alpha p_{01}^2 \ =\ \alpha p_{03}^2 
    \ =\ \beta p_{12}^2& = & \beta p_{23}^2 \, . \label{eq:eq3}
 \end{eqnarray}
For $\alpha,\beta$ satisfying~\eqref{E:abnez} as well as
$(1-\alpha)(1-\beta)\neq 0$, we distinguish the
following three disjoint cases. 
\medskip

\noindent
\emph{Case 1:} $p_{02} = 0$. 

\noindent
Since $p_{13} = 0$ would imply that all components are zero and
hence contradict $(p_{01}, \ldots, p_{23})^{\mathrm{T}} \in \P^5$, we can
assume $p_{13} = 1$. 
Then~\eqref{eq:eq1} and
(\ref{eq:eq3}) imply
\[
  \alpha p_{01}^2\ =\ \alpha p_{03}^2\ =\ 
   \beta p_{12}^2\ =\ \beta p_{23}^2\ =\ 
    \frac{\alpha\beta}{(1-\alpha)(1-\beta)} \neq 0 \, .
\]
Since (\ref{eq:eq2}) implies
$p_{01} p_{23} = -p_{03} p_{12}$,
only 8 of the $2^4=16$ sign combinations for $p_{01}, p_{03}, p_{12}, p_{23}$
are possible. 
Namely, the 8 (complex) solutions for 
$p_{01}, p_{03},p_{12}, p_{23}$ are 
\begin{equation}
\label{eq:solutions4tuple}
  (p_{01},p_{03},p_{12},p_{23})^{\mathrm{T}}
  = \frac{1}{\sqrt{(1-\alpha)(1-\beta)}} \,
    (\gamma_{01} \sqrt{\beta}, \gamma_{03} \sqrt{\beta}, 
      \gamma_{12} \sqrt{\alpha},
    - \gamma_{01} \gamma_{03} \gamma_{12} \sqrt{\alpha})^{\mathrm{T}}
\end{equation}
with $\gamma_{01}, \gamma_{03}, \gamma_{12} \in \{-1,1\}$.
Hence, for $\alpha,\beta \in \R^2$ satisfying~\eqref{E:abnez}, 
this case gives 8 distinct common tangents.

\medskip

\noindent
\emph{Case 2:} $p_{13} = 0$. \\
This case is symmetric to case 1. Setting $p_{02} = 1$, the resulting
8 solutions for the variables
$p_{01}, p_{03}, p_{12}, p_{23}$ are the same ones as
in~(\ref{eq:solutions4tuple}).
\medskip

\noindent
\emph{Case 3:} $p_{02} p_{13} \neq 0$. \\
Without loss of generality, we can assume $p_{02} = 1$.
Solving (\ref{eq:eq2}) for $p_{13}$ and substituting this expression into
$(\ref{eq:eq1})$ yields
\[
  - \beta - \beta p_{01}^2 p_{23}^2
  - \beta p_{03}^2 p _{12}^2
  - 2 \beta p_{01} p_{03} p_{12} p_{23}
  + (1-\alpha)(1-\beta) p_{03}^2\ =\ 0 \, .
\]
We use~\eqref{eq:eq3} to write this in terms of $p_{01}$.
This is straightforward for the squared terms, but for the other terms, 
we observe that, by~\eqref{eq:eq3}, $p_{01}p_{23}=\pm p_{03}p_{12}$ and since  
$p_{02} p_{13} \neq 0$, the Pl\"ucker equation~\eqref{eq:eq2} implies these
have the same sign.
This gives the quartic equation in $p_{01}$
\[
  - \beta + (1-\alpha)(1-\beta) p_{01}^2 - 4 \alpha p_{01}^4\ =\ 0 \, ,
\]
whose discriminant is
 \begin{equation}\label{E:disc}
  (1-\alpha)^2(1-\beta)^2-16\alpha\beta \, .
 \end{equation}
Hence, for $\alpha,\beta \in \R$ satisfying~\eqref{E:abnez}, and for which
this discriminant does not vanish, 
there are two different solutions for $p_{01}^2$.
For each of these two solutions for $p_{01}^2$, there are 8 distinct solutions
for $p_{01}, p_{03}, p_{12}, p_{23}$, namely
\begin{equation}
\label{eq:solutions4tuplecase3}
  (p_{01},p_{03},p_{12},p_{23})^{\mathrm{T}}
  \ =\ \sqrt{ p_{01}^2} \,
    ( \gamma_{01}, \gamma_{03}, \gamma_{12},
      \gamma_{01} \gamma_{03} \gamma_{12})^{\mathrm{T}}
\end{equation}
with $\gamma_{01}, \gamma_{03}, \gamma_{12} \in \{-1,1\}$. 
Since $p_{13}$ is uniquely determined by 
$p_{01}$, $p_{02}$, $p_{03}$, $p_{12}$,  case~3
gives 16 distinct common tangents.\medskip

In order to determine when all solutions are real,
suppose first that $\alpha=\beta$.
Then the discriminant~\eqref{E:disc} becomes
$(\alpha^2-6\alpha+1)(\alpha+1)^2$, and its smallest positive root
is $\alpha_0 := 3 - 2 \sqrt{2} \approx 0.17157$.
In particular, for $0 < \alpha < \alpha_0$, 
the discriminant in case~3 is positive and both solutions 
for $p_{01}^2$ are positive.
Thus, for $0 < \beta=\alpha < \alpha_0$,
the solutions of all three cases are distinct and real. 
Next, fix $0 < \alpha < \alpha_0$ and suppose that $0 < \beta<\alpha$.
Then the discriminant~\eqref{E:disc} is positive:
for fixed $0 < \alpha < \alpha_0$, the 
discriminant~\eqref{E:disc} is decreasing in $\beta$ for $0 < \beta < \alpha$
and positive when $\beta = \alpha$.
This concludes the proof of Theorem~\ref{th:32construction}. 
\end{proof}

Figure~\ref{fi:explproj} illustrates
the construction and the 32 tangents for $\alpha=1/10$ and $\beta=1/20$.
\ifpictures
\begin{figure}[htb]
\[
  \epsfxsize=3.3in \epsfbox{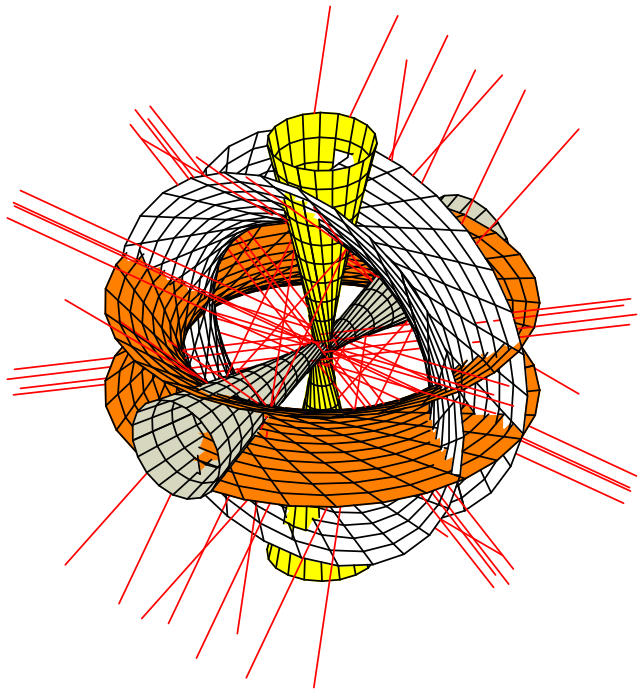}
\]
\caption{The configuration of quadrics from Theorem~\ref{th:32construction}.}
\label{fi:explproj}
\end{figure}
\fi

\subsection*{Acknowledgment.} Part of this work was done while the authors
stayed at Mathematical Sciences Research Institute, Berkeley.
 
\providecommand{\bysame}{\leavevmode\hbox to3em{\hrulefill}\thinspace} 
\providecommand{\MR}{\relax\ifhmode\unskip\space\fi MR } 
\providecommand{\MRhref}[2]{%
  \href{http://www.ams.org/mathscinet-getitem?mr=#1}{#2} 
} 
\providecommand{\href}[2]{#2}


\end{document}